\newcommand{\cc}{\mathbb{C}}

\def\EL{\mathop{\rm EL}}
\def\VEL{\mathop{\rm VEL}}
\def\EEL{\mathop{\rm EEL}}

\def\area{\mathop{\rm area}}

\def\osph{S_\mathcal{O}}
\def\osphn{\osph(n)}

\def\einc{\mathop{\rm eInc}}
\def\vinc{\mathop{\rm vInc}}
\def\einca{\mathop{\rm eIncAdj}}

\documentclass{article}

\usepackage{amsmath}
\usepackage{url}
\usepackage{graphicx}
\usepackage{latexsym}
\usepackage{amsfonts}
\usepackage{amssymb}
\usepackage[font=sf]{caption}
\usepackage{subfigure}
\usepackage{verbatim}

\newtheorem{theorem}{Theorem}[section]
\newtheorem{lemma}[theorem]{Lemma}

\newtheorem{corollary}[theorem]{Corollary}

\newcommand{\boxx}{\unskip \nopagebreak \hfill \rule{0.75em}{0.75em}}
\newenvironment{proof}%
    {\noindent{\bfseries Proof.}\begin{rmfamily}}%
    {\boxx\end{rmfamily}}

\numberwithin{equation}{section}

\begin{document}

\title{Combinatorial modulus and type of graphs}



\author{William E. Wood\\Knox College\thanks{This paper documents results of the author's Ph.D.
thesis \cite{woodthesis} prepared at Florida State University under
the direction of Philip L. Bowers, whose guidance and support were
invaluable to this work.}}
\date{October 11, 2006}

\maketitle


\begin{abstract}
Let a $A$ be the 1-skeleton of a triangulated topological annulus.
We establish bounds on the combinatorial modulus of  a refinement
$A'$, formed by attaching new vertices and edges to $A$,  that
depend only on the refinement and not on the structure of $A$
itself. This immediately applies to showing that a disk
triangulation graph may be refined without changing its
combinatorial type, provided the refinement is not too wild.  We
also explore the type problem in terms of disk growth, proving a
parabolicity condition based on a superlinear growth rate, which we
also prove optimal.   We prove  our results with no degree
restrictions in both the EEL and VEL settings and examine type
problems for more general complexes and dual graphs.
\end{abstract}

\maketitle


\section{Introduction}

There are two ways carry the notion of conformal modulus of a ring
domain to a triangulated annulus, depending on whether metrics are
assigned to the vertices or the edges of the 1-skeleton.  The two
versions are qualitatively different and even lead to inequivalent
notions of discrete conformal type  -- VEL type for vertices, EEL
type for edges.

The first goal of this paper is to establish how subdividing the
faces of a triangulated annulus can affect its discrete modulus in
either setting.  We show that the distortion of the modulus may be
bounded in terms of the subdivision alone with no dependence on the
original triangulation.  In particular, there is no dependence on
degree. This is avoided by applying an observation of Chrobak and
Eppstein \cite{eppstein} that  any planar graph may be considered as
a directed graph with globally bounded outdegree. This weaker notion
of bounded degree is sufficient to control the bounds.

 This has immediate application
to discrete type problems, i.e. determining whether a disk
triangulation graph is hyperbolic or parabolic in either the VEL or
EEL setting. We show that if a disk triangulation graph is refined
in a sufficiently reasonable way, the resulting graph will have the
same type. We obtain different notions of ``sufficiently
reasonable'' for VEL and EEL types, but both will cover most
standard refinement processes, such as hexagonal and barycentric
subdivision.

We then turn to the exploring discrete type in terms of the growth
of spheres.  There are already some results of this ilk (e.g.,
\cite{rs}, \cite{siders}, \cite{soardi}), but they require symmetry
or degree restrictions on the graph.  We obtain a superlinear growth
condition that guarantees parabolicity with no such restrictions.
We also show how to construct slow-growing (e.g.,
$n^{1+\varepsilon}$) hyperbolic graphs, establishing sharpness of
the parabolicity condition.

We establish our definitions and foundational lemmas in
Section~\ref{sec:prelim}.  Our main results regarding the moduli of
ring domains are developed and proved in Section~\ref{sec:refine}.
We apply these results to the type problem in Section~\ref{sec:type}
and offer examples demonstrating the necessity of our hypotheses. We
show how to generalize our results to non-triangular complexes in
Section~\ref{sec:dual} and apply this result to relate the type of a
complex to that of its dual.  We also introduce discrete outer
spheres and explore their application to discrete type problems. Our
results relating type to sphere growth are covered in
Section~\ref{sec:growth}.

\section{Preliminaries}\label{sec:prelim}

\subsection{Extremal length} Our definitions for combinatorial extremal length
are  consistent with \cite{heschrammhp}.

Let $X$ be a non-empty set and $\Gamma$ a non-empty collection of
finite or infinite sequences in $X$, called  \emph{paths}.  We will
thus refer to the \emph{pair} $(X,\Gamma)$.  A \emph{metric} on $X$
is a function $m:X\to [0,\infty)$. The value $m(x)$ is the
\emph{$m$-weight} or \emph{$m$-measure} of $x$. The \emph{area} of
$m$ is
$$\mathrm{area}(m) = \sum_{x\in X} m(x)^2,$$ and a metric is called \emph{admissible} if it has finite, non-zero area.  Let
${\mathcal{M}(X)}=\{m:\area(m) <\infty\}$ be the set of admissible
metrics on $X$.  For a path $A=\{a_0,a_1,\ldots\}\subset X$, define
its $m$-\emph{length} to be $L_m(A) = \sum_{j=1}^\infty m(a_j)$.  We
abuse notation by writing for convenience $\sum_{x\in A}
m(x)=\sum_{j=1}^\infty m(a_j)$. For a collection $\Gamma$ of paths
in $X$, define $L_m(\Gamma) = \inf_{A\in \Gamma} L_m(A)$ and the
extremal length
$$\EL(\Gamma) = \sup_{m\in \mathcal{M}(X)} \left\{
\frac{L_m(\Gamma)^2}{\area(m)}\right\}.$$ The reciprocal of extremal
length is the \emph{modulus}.

We say $\Gamma$ is \emph{hyperbolic} if $\EL(\Gamma)$ is finite and
\emph{parabolic} if $\EL(\Gamma)$ is infinite.  When the set
$\Gamma$ is clear from the context, we  refer to
$\EL(X)=\EL(\Gamma)$.

An \emph{extremal metric} for $\Gamma$ is an admissible metric $\mu$
on $X$ for which $\EL(\Gamma)=\frac{L_\mu(\Gamma)^2}{\area(\mu)}$.
 In the case $\EL(\Gamma)=\infty$, an
extremal metric has finite area and all elements of $\Gamma$ have
infinite length.
\begin{lemma} Let $X$ be set and $\Gamma$ a collection of subsets.  If  $\Gamma$ is finite or if $EL(\Gamma)=\infty$, then there is an extremal metric
for $\Gamma$ on $X$. \end{lemma} The finite case is proved in
\cite{cannonacta}.  The latter case defines a \emph{parabolic
extremal metric}; its existence is an exercise in
\cite{heschrammhp}.

Note that scaling the metric does not change the quantity maximized
by extremal length and so we may assume that our metrics are always
normalized to have area one.

He and Schramm also offer in \cite{heschrammhp} the important
\emph{monotonicity property}, stated as
\begin{lemma}\label{th:monotonicity}  Suppose $\Gamma $ and $\Gamma'$ are collections of subsets of $X$ with the
property that for every $\gamma\in\Gamma$ there is a
$\gamma'\in\Gamma'$ such that $\gamma'\subset\gamma$. (In
particular, this holds if $\Gamma\subset\Gamma'$.) Then
$\EL(\Gamma')\leq \EL(\Gamma)$.
\end{lemma}

\subsection{Comparability}

Let $f$ and $g$ be positive real-valued functions with domain
$\Upsilon$.  Let $k\geq 1$. We say the functions are
\emph{$k$-comparable} if $\frac{1}{k}g(\upsilon)\leq f(\upsilon)\leq
kg(\upsilon)$ for all $\upsilon\in\Upsilon$. $f$ and $g$ are
\emph{comparable} if they are $k$-comparable for some $k\geq 1$. It
is easy to verify that comparability defines an equivalence
relation, and this is the relation we  seek when determining
finiteness of extremal length.  Its application is prescribed by the
following lemma.

\begin{lemma} Let $X$ and $Y$ be infinite sets with sets of paths $\Gamma_X$ and $\Gamma_Y$. Let
$\{A_i\}_{i=0}^\infty$ and $\{B_i\}_{i=0}^\infty$ be collections of
finite subsets of $X$ and $Y$, respectively, and $\Gamma^i_X=\{
\gamma\cap A_i: \gamma\in \Gamma_X\}$ and $\Gamma^i_Y=\{ \gamma\cap
B_i: \gamma\in \Gamma_Y\}$.  Suppose these sets satisfy the
following properties:\begin{enumerate} \item $X=\cup_{i\geq 0} A_i$
and $Y=\cup_{i\geq 0} B_i.$ \item If $i<j$, then $A_i\subset A_j$
and $B_i \subset B_j$. \item Let $i\geq 0$. For every $\gamma_X\in
\Gamma_X$ and $\gamma_Y\in \Gamma_Y$, $\gamma_X\cap A_i$ and
$\gamma_Y\cap B_i$ are non-empty. \item $\EL(\Gamma_X^i)$ and
$\EL(\Gamma_Y^i)$ are comparable (taken as functions of
$i$).\end{enumerate} Then $\EL(\Gamma_X)=\infty$ if and only if
$\EL(\Gamma_Y)=\infty$.\label{th:comparabletype}
\end{lemma}
\begin{proof}
Suppose $\EL(\Gamma_X)=\infty$ with parabolic extremal metric $\mu$.
Define $\mu_i$ to be the restriction of $\mu$ to $A_i$ and note that
$\area(\mu_i)\leq \area(\mu)=1$.
 Choose any $N>0$. We show $\EL(\Gamma_Y)>N$, implying
$\EL(\Gamma_Y)=\infty$.

Since $X$ is parabolic, every element $\gamma\in \Gamma_X$ has
infinite $\mu$-length.  Choose $k>0$ so that $k\EL(\Gamma_X^i)\leq
\EL(\Gamma_Y^i)$ for all $i>0$.  All paths in $\Gamma_X$ have
infinite $\mu$-length, and so for any given path $\gamma\in\Gamma_X$
there is a $j_\gamma>0$ so that $L_{\mu_{j_\gamma}}(\gamma\cap
A_j)>\sqrt{\frac{N}{k}}$.  Let $j=\inf_{\gamma\in\Gamma_X}
j_\gamma$. Since every path in $\Gamma_X^j$ is contained in a
transient path in $\Gamma_X$, we have $L_{\mu_j}(\Gamma^j_X)>
\sqrt{\frac{N}{k}}$. Then
$$N< k L_{\mu_j}(\Gamma^j_X)^2 \leq
k\frac{L_{\mu_j}(\Gamma^j_X)^2}{\area(\mu_j)}$$
$$ \leq
k\sup_{m\in\mathcal{M}(A_j)}\frac{L_{m}(\Gamma^j_X)^2}{\area(m)}=k\EL(\Gamma_X^j)\leq
\EL(\Gamma_Y^j)\leq \EL(Y).$$ The last inequality is a direct
consequence of monotonicity (Lemma~\ref{th:monotonicity}). The proof
is completed by repeating the argument with the roles of $X$ and $Y$
exchanged.

\end{proof}

\section{Refinement and extremal length}\label{sec:refine}

\subsection{Shadow paths}\label{sec:shadowpaths}

Our goal is to control  the combinatorial extremal length of a set
$X$ that is related to some other set $X'$ whose extremal length is
known. We codify our technique in the following lemma.

\begin{lemma} Let $X'$ be a set, $\Gamma'$ a collection of subsets
of $X'$, and $\mu'$ an extremal metric for $(X',\Gamma')$.  Suppose
there is a set $X$, a collection $\Gamma$ of subsets of $X$, an
admissible metric $\mu$ on $X$, and constants $C,D>0$ with the
following properties:
\begin{enumerate}
\item $\area(\mu)\leq C\cdot \area(\mu')$
\item For each $\gamma\in\Gamma$, there is a $\gamma'\in\Gamma'$
such that $D\cdot L_{\mu'}(\gamma')\leq
L_\mu(\gamma)$.\end{enumerate} Then
$\EL(\Gamma)\geq\frac{D^2}{C}\EL(\Gamma').$\label{th:shadowpaths}
\end{lemma}
\begin{proof} We associate to each
$\gamma\in\Gamma$ a specific path $\gamma'\in\Gamma'$ with
$L_{\mu}(\gamma)\geq D\cdot L_{\mu'}(\gamma')$. Let $\Gamma^\#$ be
the collection of these $\gamma'$.  The proof now amounts to
unraveling the definitions.
$$\EL(\Gamma)=\sup_{m\in\mathcal{M}(X')}\frac{\inf_{\gamma\in\Gamma}
L_m(\gamma)^2}{\area(m)}\geq \frac{\inf_{\gamma\in\Gamma}
L_\mu(\gamma)^2}{\area(\mu)}$$
$$\geq \frac{\inf_{\gamma'\in\Gamma^\#}
(D\cdot L_{\mu'}(\gamma'))^2}{\area(\mu)} \geq D^2
\frac{\inf_{\gamma'\in\Gamma'} L_{\mu'}(\gamma')^2}{\area(\mu)}$$
$$\geq D^2\frac{\inf_{\gamma'\in\Gamma'}
 L_{\mu'}(\gamma')^2}{C\cdot\area(\mu')}
=\frac{D^2}{C}\cdot\frac{\inf_{\gamma'\in\Gamma'}
L_{\mu'}(\gamma')^2}{\area(\mu')}=\frac{D^2}{C}\EL(\Gamma').$$
\end{proof}

The set $\Gamma^\#$ is the set of \emph{shadow paths} and is
essential to the forthcoming results. Suppose we want to find the
extremal length of a pair $(X,\Gamma)$ that is constructed from
another pair $(X',\Gamma')$ whose extremal length is known. Assume
an extremal metric $\mu'$ on $X'$.  We  then use $\mu'$ to construct
a new metric $\mu$ on $X$ satisfying the assumptions of
Lemma~\ref{th:shadowpaths}, meaning we always have two things to
control: area and path length.  The trick is to construct the metric
so that the constants $C$ and $D$ depend on as little as possible.

For $x\in X$, $\mu(x)$ is assigned a value of $\mu(x')$ for some
$x'\in X'$.  That is, each element $x\in X$ has a corresponding
element $x'\in X'$ that prescribes its measure. The constant $C$ is
a bound on the number of elements in $X$ to which an element of $X'$
may be assigned.

For a path $\gamma$ in $\Gamma$, we must guarantee a path in
$\gamma'\in\Gamma'$ whose $\mu'$-length is less than $\frac{1}{D}$
times the $\mu$-length of $\gamma$.  The paths $\gamma$ and
$\gamma'$ naturally correspond.  As $\gamma$ bobs and weaves through
$X$, $\gamma'$ will ``shadow'' its movement in $X'$ and have
comparable length.

These two conditions are at odds.  We need to choose $\mu$ carefully
so that paths are sufficiently long, but so that the area stays
sufficiently small.

Our objective is comparability of the extremal lengths of two sets
$A$ and $B$.  This requires applying Lemma~\ref{th:shadowpaths}
twice, with $A$ and $B$ alternatively taking the roles of $X$ and
$X'$. The extremal lengths of $A$ and $B$ are shown to be
$k$-comparable for some $k$ depending only on the constants in the
lemma.

\subsection{Graphs}

A \emph{graph} $G=(V,E)$ is a set $V$ of vertices and a set $E$ of
edges.  An edge connecting two vertices $v_0$ and $v_1$ is denoted
$[v_0,v_1]$, and the vertices $v_0$ and $v_1$ are the
\emph{endpoints}. A graph is finite or infinite as $|V|$ is finite
or infinite. Note that the definition and notation disallow multiple
edges connecting two vertices.  We assume our graphs have no edges
connecting a vertex to itself, and all of our graphs are connected.
Two vertices $v,w$ are \emph{adjacent}, denoted $v\sim w$, if
$[v,w]\in E$. Two graphs $G=(V,E)$ and $G'=(V',E')$ are
\emph{isomorphic}, denoted $G\cong G'$, if there is a bijection
$f:V\to V'$ such that for any $v,w\in V$ we have $v\sim w$ if and
only if $f(v)\sim f(w)$. Every vertex has a \emph{degree}
$\mathrm{deg}(v_0) = |\{v\in V: [v,v_0]\in E\}|$, the number of
edges having $v_0$ as an endpoint. We will assume the degree is
finite for each vertex (the graph is \emph{locally finite}).  The
degree of a graph is defined as $\sup_{v\in V}\mathrm{deg}(v)$.  The
graph has \emph{bounded degree} if the degree is finite,
\emph{unbounded degree} otherwise.

A \emph{vertex path} in $G$ is a finite or infinite sequence of
vertices $v_0, v_1, \ldots$ such that for all $i\geq 0$, either
$v_i\sim v_{i+1}$ or $v_i = v_{i+1}$. Similarly, an \emph{edge path}
is a sequence of edges of the form $[v_0, v_1],[v_1, v_2],[v_2,
v_3],\ldots$, i.e. consecutive edges laid end to end. Note that each
vertex path has a corresponding edge path, and vice versa. The
combinatorial vertex length and edge length of a path are the
numbers of vertices and edges in the path, respectively.  If a base
point $v_0\in V$ is specified, we say the \emph{norm} of a vertex
$w\in V$ is the combinatorial vertex length of the shortest path
connecting $v_0$ to $w$. A \emph{cycle graph} is a finite connected
graph whose vertices all have degree two.

We generally think of the vertices as points and the edges as arcs
connecting the endpoints.   A graph is \emph{planar} if the graph
can be embedded in the plane -- that is, the vertex points and edge
arcs may be positioned in the plane so that the vertices are located
at distinct points and the edge arcs intersect only at shared
endpoints. An embedding like this is a \emph{diagram} for a graph.
For example, the diagram of a cycle is a Jordan curve.

Our interest in graphs is to mimic the geometric properties of
classical Riemann surfaces with a combinatorial object.  This is
achieved via the \emph{triangulation graph}, which is  the
1-skeleton of a locally finite tiling by triangles of a simply
connected Riemann surface (possibly with boundary). A triangle with
vertices or edges $a,b,c$ is denoted $\triangle(a,b,c)$. Of specific
interest are \emph{disk triangulation graphs}, in which the Riemann
surface in question is the open unit disk.  We  also study more
general \emph{disk cell complexes} whose 1-skeletons are \emph{disk
cell graphs}. These are obtained from locally finite tilings of the
disk by general finite polygons, not necessarily triangles.  The
interiors of the polygons are the \emph{faces} or \emph{cells}. We
frequently require a global bound on the number of sides of the
polygonal faces. When a cell structure is present and relevant, we
will expand the graph notation $G=(V,E,F)$ to include the set of
polygonal faces $F$.

A \emph{directed graph} is a graph $G=(V,E)$ along with an ordering
on the vertices of each edge.  For an edge $[x,y]\in E$, we may form
a \emph{directed edge} by $[x,y\rangle$, where $x$ is the \emph{tail
vertex} and $y$ is the \emph{head}.  For any vertex $v\in V$, the
\emph{outdegree} of $v$ is the number of directed edges for which
$v$ is the tail vertex.

Let $G$ be a cell complex embedded in the plane with disjoint
subcomplexes $A,B,C\subset G$. We say $C$ \emph{separates} $A$ from
$B$ if $A$ and $B$ lie in different components of $G\setminus C$. We
say $C$ \emph{separates $A$ from infinity} if $A$ lies in an
unbounded component of $G\setminus C$.

 We say a
planar graph is a \emph{combinatorial annulus} or an \emph{annular
complex} if the graph is the 1-skeleton of triangulated topological
annulus in the plane, i.e. a region $A\subset\cc$ whose boundary
consists of two disjoint Jordan curves $C_1,C_2$ such that the
bounded component of $\cc\setminus C_2$ contains the bounded
component of $\cc\setminus C_1$.

We may apply the definition of combinatorial extremal length to a
graph $G=(V,E)$ in two ways, according to whether the metric assigns
values to the set of vertices or to the set of edges (\emph{vertex
metrics} and \emph{edge metrics}).

Let $G$ be a combinatorial annulus and let $\Gamma_V=\Gamma_V(G)$ be
the set of vertex paths connecting the two boundary components of
$G$.  Define $\Gamma_E=\Gamma_E(G)$ similarly for edge paths. Define
the \emph{vertex extremal length} $\VEL(G) = \EL(\Gamma_V)$, and the
\emph{edge extremal length} $\EEL(G)=\EL(\Gamma_E)$. A graph is
\emph{VEL-hyperbolic} or \emph{VEL-parabolic} as the vertex extremal
length is finite or infinite, indicating its \emph{VEL type}.
\emph{EEL-hyperbolic}, \emph{EEL-parabolic}, and \emph{EEL type} are
defined similarly.

\subsection{Refinement}\label{sec:refinement}

Let $G=(V,E,F)$, $rG=(rV,rE,rF)$ be planar complexes and suppose
there is an injection $\iota:V\to rV$ with the property that if
$[x,y]\in E$, there is a collection of vertices
$w_0=\iota(x),w_1,\dots,w_n=\iota(y)$ with $[w_j,w_{j+1}]\in rE$ and
$w_j \notin \iota(V)$ for each $j\neq 0,n$. We then call $rG$ a
\emph{refinement} of $G$ and we shall consider the refinement $r$ as
a map from an appropriate set of  graphs to itself such that $rG$ is
always a refinement of $G$. We abuse notation by suppressing further
mention of the injection $\iota$ and considering $V\subset rV$.

Less technically, a refinement attaches vertices to the edges of a
graph, dividing the edge into subedges, and then adds edges inside
the faces.

\begin{figure}\centering
    \includegraphics[width=1in]{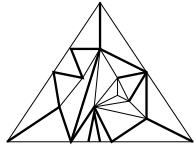}
 \caption{Bold edges are incident-adjacent.}\label{fig:incident}
\end{figure}

We  require some notation before specifying the classes of
refinements we need to consider.
 Let $G=(V,E,F)$ and $r$ a refinement of $G$. For $e=[x,y]\in E$, we
refer to the vertices in $rV$ \emph{incident} to $e$ as the set
$\vinc_r(e)=\{w_1,\dots,w_{n-1}\}$ guaranteed by the definition of
refinement (note that we do not include the endpoints), and
similarly the set of edges incident to $e$ is
$\einc_r(e)=\{[w_0,w_1],\dots,[w_{n-1},w_n]\}$. We  also consider
edges that are \emph{incident-adjacent}, $\einca_r(e)$
$=\{[x',y']\in rE:$ $ x'\in\vinc_r(e)\cup\{x,y\}$ and $[x',y']$ is
contained in a cell bounded by $e\}$.  This definition requires $G$
to have a planar cell structure. See Figure~\ref{fig:incident}.  For
any vertex $v\in rV\setminus V$ not incident to an edge of $G$,
there is a face $F$ of $G$ whose incident edges separate $v$ from
infinity. We say $v$ then \emph{lies in} $F$.  Similarly, an edge
lies in $F$ if either of its endpoints lies in $F$ and it is not
incident to any edge of $G$.

\begin{figure}\centering\label{fig:refinementexamples}
 \subfigure[]{
     \includegraphics[width=1in]{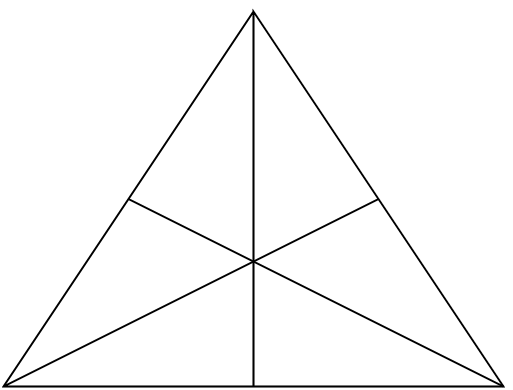}
        }
 \subfigure[]{
        \includegraphics[width=1in]{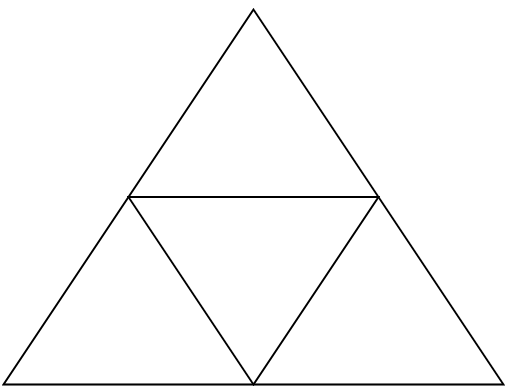}
        }
 \caption{(a) Barycentric and (b) hexagonal refinement}
\end{figure}

We  examine two families of refinement.  A refinement $r$ is
$b$-\emph{weakly bounded} if $|\einc_r(e)|\leq b$ for all $e\in E$,
and $r$ is $(b,c)$-\emph{strongly bounded} if it is $b$-weakly
bounded and  for any edge $e\in E$ and any vertex $v'\in V'$
bounding $e$ or incident to $e$ under $r$, we have that
$\einca_r(e)$ contains at most $c$ edges attached to $v'$ lying in
any one cell. Basically, weakly bounded means that $r$ partitions
the edges of $G$ into at most $b$ pieces. Strongly bounded further
assumes a bound on the number of new edges attached to the boundary
of each cell of $G$.  Note that if $G$ does not have bounded degree,
then it is possible in a  strongly bounded refinement for the number
of edges attached to a vertex to be unbounded; it is bounded within
each cell, but there may be unboundedly many cells.  Examples of
strongly bounded refinements include the identity refinement,
barycentric subdivision, and hexagonal refinement.  See
Figure~\ref{fig:refinementexamples}.

These are our notions of ``sufficiently reasonable'' mentioned in
the introduction.  Weakly and strongly bounded refinements will
naturally correspond to vertex and edge extremal length,
respectively. Indeed, the method He and Schramm offer in
\cite{heschrammhp} to construct a VEL-parabolic, EEL-hyperbolic
triangulation amounts to a strongly bounded refinement of a
VEL-parabolic graph. VEL type is the same as circle packing type,
whereas EEL type indicates recurrence or transience of a random walk
or electric network (see \cite{duffin}, \cite{doyle},
\cite{heschrammhp}, \cite{ken}).

\subsection{Edge extremal length}

\begin{theorem} \label{th:eel} Let $A$ be a finite combinatorial annulus and $r$ a $(b,c)$-strongly bounded
refinement. Then there is a constant $k\geq 1$, depending only on
$b$ and $c$, such that $\EEL(rA)$ and $\EEL(A)$ are $k$-comparable.
That is,
$$\frac{1}{k} \EEL(rA) \leq \EEL(A) \leq k \EEL(rA).$$
\end{theorem}

\begin{proof}
Suppose $\mu'$ is an edge extremal metric on $rA$. We show the
function
$$\mu(e)=\max_{e'\in \einc_r(e)} \mu'(e')$$  defines an edge metric on $A$ that provides the necessary bound.
First,
$$\area(\mu) =\sum_{e\in E} \mu(e)^2 = \sum_{e\in E} \max_{e'\in \einc_r(e)}
\mu'(e')^2 \leq \sum_{e'\in rE} \mu'(e')^2 =\area(\mu').$$  The
inequality holds because no edge in $rE$ is incident to more than
one edge in $E$.

Now suppose $\gamma\in\Gamma=\Gamma_E(A)$.  Construct the shadow
path $\gamma'\in\Gamma'=\Gamma_E(rA)$ replacing each $e\in\gamma$
with the appropriately ordered edges of $\einc_r(e)$. Then
$$L_\mu(\gamma)=\sum_{e\in\gamma}\mu(e) =\sum_{e\in\gamma} \max_{e'\in
\einc_r(e)} \mu'(e')  =\frac{1}{b}\sum_{e\in\gamma} b \max_{e'\in
\einc_r(e)} \mu'(e') $$ $$\geq \frac{1}{b}\sum_{e\in\gamma}
\sum_{e'\in
 \einc_r(e)} \mu'(e) = \frac{1}{b} \sum_{e'\in\gamma'}\mu'(e)= \frac{1}{b}
L_{\mu'}(\gamma')$$ and so by Lemma~\ref{th:shadowpaths}
$$\EEL(A) \geq \frac{1}{b^2}\EEL(rA).$$

Conversely, suppose $\mu$ is an edge extremal metric on $A$ and
define an edge metric $\mu'$ on $rA$ by
\begin{displaymath}\mu'(e') = \left\{
  \begin{array}{ll}
     \mu(e) & \textrm{if there is an $e\in E$} \\ & \textrm{such that $e'\in  \einc_r(e)$.}
     \\ \\
     \max(\mu(e_1), \mu(e_2), \mu(e_3)) & \textrm{if } e'\subset\triangle(e_1,e_2,e_3) \textrm{ and} \\
     & e' \textrm{ is incident-adjacent} \\
     & \textrm{to one of the $e_1,e_2,$ or $e_3$.} \\ \\
     0 & \textrm{otherwise.}
  \end{array}
  \right.
\end{displaymath}
Then
$$\area(\mu')=\sum_{e'\in rE} \mu'(e')^2 \leq \sum_{e\in
E}\sum_{e'\in\einc_r(e)} \mu(e)^2+\sum_{e\in E}6c \mu(e)^2 $$ $$\leq
b \area(\mu) +6(1+b)c\area(\mu)=\left(b+6c(1+b)\right)\area(\mu).$$
The first term of the inequality comes from the fact that every edge
in $E$ contributes its measure to at most $b$ incident edges.  The
second term counts the incident-adjacent edges. An edge $e\in E$
lies on the boundary of two cells $\tau_1, \tau_2$ and may
contribute its measure to elements of $rE$ that are
incident-adjacent to any of the six edges  bounding $\tau_1$ and
$\tau_2$ (double counting $e$, once for each cell it bounds) and are
contained within one of these cells. That makes six edges, each
partitioned at most $b$ times, and at most $c$ incident-adjacent
edges attached to any refined vertex lying in a given cell.  Noting
that $b$ only counts the new vertices in the refinement, we also get
another possible $6c$ incident-adjacent edges off of the original
vertices of the $\tau_i$, making a total contribution of $6c(1+b)$.

Now suppose $\gamma'\in\Gamma'$.  Write $\gamma'$ as a concatenation
of segments $\gamma'_0,\sigma'_0,$ $\gamma'_1,\sigma'_1\dots$
$\gamma'_n, \sigma'_n\dots$ such that for every $i\geq 0$ there is
an $e_i\in E$ with $\sigma_i\subset\einc_r(e_i)$ and there exist
$e^i_1,e^i_2\in E$ (not necessarily distinct) such that the two end
edges  of $\gamma'_i$ are contained in $\einca_r(e^i_1)$ and
$\einca_r(e^i_2)$, and all edges in $\gamma'_i$ lie inside the same
triangular cell. (Some of the $\sigma_i$ and $\gamma_i$ may be
empty.)  Define $\gamma_i$ to be the edge path $\{e^i_1, e^i_2\}$
and $\sigma_i=\{e_i\}$. Let $\gamma$ be the concatenation
$\gamma_0\sigma_0\gamma_1\sigma_1\ldots$. Then
$$L_{\mu'}(\gamma')=\sum_{i\geq 0} L_{\mu'}(\gamma'_i)+\sum_{i\geq 0} L_{\mu'}(\sigma'_i)=\sum_{i\geq 0}
\big(\sum_{x\in\gamma'_i} \mu'(x)+ \sum_{y\in\sigma'_i}
\mu'(y)\big)$$
$$\geq \sum_{i\geq 0} (\mu(e^i_1)+\mu(e^i_2))+\sum_{i\geq 0}\mu(e_i)$$ $$=\sum_{i\geq 0}L_\mu(\gamma_i)+\sum_{i\geq 0}
L_\mu(\sigma_i)= L_\mu(\gamma).$$ We now apply
Lemma~\ref{th:shadowpaths} to obtain
$$\EEL(A) \geq \frac{1}{b+6c(1+b)}\EEL(rA).$$
The proof is completed by taking $k=\max(b^2, b+6c(1+b))$.
\end{proof}

We have actually proved a little more than is stated, since the
constant in first half of the proof did not depend on $c$.  As such,
the following is a corollary of the proof.
\begin{lemma} Let $A$ be a finite combinatorial annulus and $r$ a $b$-weakly bounded
refinement. Then $ \EEL(rA) \leq b^2 \EEL(A).$
\end{lemma}

The point of this and subsequent theorems is that the comparability
constant does not depend at all on $A$. Only  the process by which
$A$ is refined matters, and even that dependence is surprisingly
slight.  This becomes important when we apply the theorem to the
type problem.

\subsection{Vertex extremal length}

We would like to take the same strategy for vertex extremal length
as we did for edge extremal length: for two annuli related by a
bounded refinement, use an extremal metric on one to construct a new
metric on the other that adequately controls its extremal length.
Unfortunately in the VEL case, the degree of the original graph
naturally arises in the bounds.  To prove a result that is not
dependent on degree, we will use the following degree property
shared by all planar graphs that will be sufficient to construct the
required metrics.

\begin{lemma} \label{th:outdeg3} The edges of a finite planar graph $G=(V,E)$ may be
directed so that the outdegree of every vertex is at most three.
\end{lemma}
See \cite{eppstein} for a proof.

\begin{theorem} \label{th:vel} Let $A$ be a finite combinatorial annulus and $r$ a $b$-weakly bounded
refinement. Then there is a constant $k\geq 1$, depending only on
$b$, such that $\VEL(rA)$ and $\VEL(A)$ are $k$-comparable.
 That is, $$\frac{1}{k} \VEL(rA) \leq \VEL(A) \leq k \VEL(rA).$$
\end{theorem}
\begin{proof}
Let $\mu'$ be an extremal metric on $rA.$  Define  a vertex metric
$\mu(v)$ on $V$ to be the larger of $\mu'(v)$ and the maximal value
of $\mu'(v')$ taken over all vertices $v'\in V'\setminus V$ that are
incident to an edge containing $v$. For any $v'\in rV$, this process
assigns $\mu'(v')$ to at most two vertices in $V$ (the two vertices
bounding the edge on which $v'$ lies) and so clearly $\area(\mu)
\leq 2\area(\mu')$.

Let $\gamma\in\Gamma(A)$ and consider the shadow path
$\gamma'\in\Gamma(rA)$ obtained by traveling between vertices of
$\gamma$ along the original edges of $A$, passing through the
incident vertices. Write $E_\gamma$ for the set of edges in $A$
connecting the vertices of $\gamma$. Then
$$L_{\mu'}(\gamma')=\sum_{v'\in\gamma'}\mu'(v')= \sum_{e\in E_\gamma}
\sum_{v'\in\vinc(e)} \mu'(v')+\sum_{v\in\gamma'\cap V} \mu'(v)$$
$$\leq \sum_{e\in E_\gamma} b \max_{v'\in\vinc(e)}\mu'(v') +\sum_{v\in\gamma'\cap V} \mu(v)\leq
b\sum_{v\in\gamma}\mu(v)+\sum_{v\in\gamma}\mu(v)=(b+1)
L_\mu(\gamma).$$ We have thus satisfied the hypotheses of
Lemma~\ref{th:shadowpaths} and conclude $\VEL(A)\geq
\frac{1}{2(b+1)^2} \VEL(rA)$.

Our work with outdegree bounds pays off in the converse. Let $\mu$
be an extremal metric on $A$ and direct the edges of $G$ so that the
outdegree of every vertex is at most 3.  This is possible by Lemma
\ref{th:outdeg3}. Define a metric $\mu'$ on $rG$ by

\begin{displaymath}\mu'(v) = \left\{
  \begin{array}{ll}
     \mu(v) & \textrm{if $v\in V$} \\ \\
     \mu(w) & \textrm{if $v\in rV\setminus V$ and $v$ is incident to} \\
     & \textrm{a directed edge with tail $w$.}\\ \\
     0 & \textrm{otherwise.}
  \end{array}
  \right.
\end{displaymath}
The bounded outdegree gives us our area restriction:
$$
\area(\mu') = \sum_{v\in rV} \mu'(v)^2 = \sum_{v\in V} \mu'(v)^2 +
\sum_{v\in rV\setminus V} \mu'(v)^2 $$ $$\leq \sum_{v\in V} \mu(v)^2
+ \sum_{v\in V} 3b \mu(v)^2
 =(1+3b) \area(\mu)=1+3b.$$

Let $\gamma'\in\Gamma(rA)$.  We  need to find  a shadow path
$\gamma\in\Gamma(A)$ with  $k L_{\mu}(\gamma) \geq
L_{\mu'}(\gamma')$ for a $k$ depending only on $b$.

Construct $\gamma$ inductively.  Suppose $\gamma'_i$ is an initial
segment of $\gamma'$ and assume $\gamma_i$ has been constructed so
that $L_{\mu'}(\gamma'_i)= L_{\mu}(\gamma_i)$.  Let $v$ be the
endpoint of $\gamma_i$ and $v'$ the endpoint of $\gamma'_i$. If
$v'\in V$, assume $v' = v$. Otherwise, $v'\in rV\setminus V$ and we
assume that if $v'$ is incident to  an edge, then that edge has
 $v$ as an endpoint. We want to extend $\gamma_i$ by a  vertex $w$ to
create a new segment $\gamma_{i+1}$ so that these properties are
preserved and so that $L_{\mu'}(\gamma'_i)=L_\mu(\gamma_i)$. The
basic rule of the construction is \emph{``always move to the tail of
the arrow.''} That is, we add to $\gamma$ the vertex at the tail of
the directed edge every time $\gamma'$ hits an incident vertex.
Figure~\ref{fig:shadowpath} illustrates the process.

Let $w'$ be the endpoint of $\gamma_{i+1}$.  There are three cases.
If $w'$ is not incident to any edge, then $\mu'(w')=0$ and we do
nothing. Set $\gamma_{i+1}=\gamma_i$ and note
$L_{\mu}(\gamma_{i+1})=L_{\mu'}(\gamma'_{i+1})$ because we do not
add any $\mu$-length.

If $w'\in V$, we take $w=w'$.  This is legal by the assumption that
$v'$ and $v$ are in the same face and that $G$ is a triangulation
(this is what guarantees $v\sim w$). Then
$L_\mu(\gamma_{i+1})=L_{\mu'}(\gamma'_{i+1})$ because we are adding
the same measure to both.

The final case is that $w'$ is incident to an edge $e$ of a face
containing $v$ as vertex.  In this case, we take $w$ to be the tail
of the directed edge $e$ ($\gamma$ may move or ``sit and wait'' at
some vertex). Again, $L_\mu(\gamma_{i+1})=L_{\mu'}(\gamma'_{i+1})$
because we are adding the same measure to both paths.

\begin{figure}\centering
\includegraphics[width=2in]{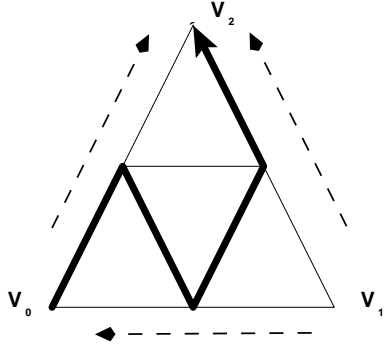}\label{fig:shadowpath}
\caption[Shadow paths in hex refinement]{An example of the shadow
path in hex refinement.  The bold arrow indicates five vertices in a
path through a hex refined face.  The dashed arrows identify how the
edges are directed.  The corresponding shadow path is $v_0, v_0,
v_1, v_1, v_2$.  Both path segments have the same length in their
respective metrics.}
\end{figure}

Thus we have constructed $\gamma$ so that $L_\mu(\gamma)=
L_{\mu'}(\gamma').$  Lemma~\ref{th:shadowpaths} now applies to give
$\VEL(rA)\geq \frac{1}{1+3b} \VEL(A)$. This proves the theorem with
$k=\max(1+3b,2b^2)$.
\end{proof}

\section{Refinement and type}\label{sec:type}

\subsection{Bounded refinement preserves type}

Our main application is to the type problem.

\begin{theorem}\label{th:type}  Let $G$ be a disk triangulation graph. If $r$ is a weakly bounded refinement, then $G$
and $rG$ have the same VEL type. If $r$ is a strongly bounded
refinement, then $G$ and $rG$ have the same EEL type.
\end{theorem}
\begin{proof} We state the proof for the VEL case.  The EEL case is identical.

  Let $C$ be the vertex cycle formed from the neighbors of the base vertex $v_0$.
Let $\{A_i\}_{i=0}^\infty$ be a collection of combinatorial annuli,
each with innermost boundary component $C$, such that $G\setminus
\cup_{i\geq 0} A_i = \{v_0\}$.  Apply Theorem~\ref{th:vel}
(Theorem~\ref{th:eel} for EEL) to conclude that $\VEL(A_i)$ and
$\VEL(rA_i)$ are comparable. The collections $\{A_i\}_{i=0}^\infty$
and $\{rA_i\}_{i=0}^\infty$  satisfy the hypotheses of
Lemma~\ref{th:comparabletype}, which says $\cup_{i\geq 0} A_i$ and
$\cup_{i\geq 0} rA_i$ have the same type. The theorem
follows.\end{proof}

\subsection{Unbounded refinements} We now present some examples illustrating the necessity of bounded
 refinements in preserving type.

\begin{theorem}  \label{th:hyprefinece} Every disk triangulation graph $G$ has a hyperbolic refinement $\zeta
G$.\end{theorem}

\begin{figure}%
\begin{center}
        \includegraphics[width=3in]{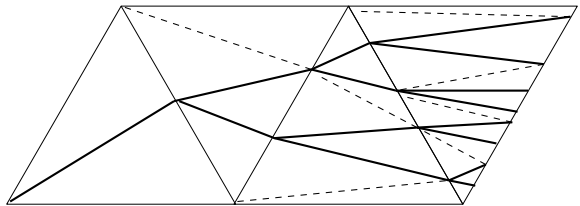}
\end{center}
\caption{Adding a binary tree with an unbounded refinement.
 The first four triangles of the refinement are shown.  Bold edges show the binary tree.
 The dashed edges are added to make the refined graph a triangulation.\label{fig:hyprefine}}
\end{figure}

\begin{proof}
Let $G$ be a VEL- or EEL-parabolic disk triangulation graph and let
$T_1, T_2,\dots$ be an infinite collection of distinct faces such
that for each $i>0$, $T_i$ intersects $T_{i+1}$  along a single edge
$e_{i+1}$, and so that $T_i$ and $T_j$ share an edge only if $|i-j|=
1$. Let $v_0$ be the vertex of $T_1$ that is not in $e_1$.

Refine the $T_i$ by attaching $2^i$ vertices to each $e_i$ and
connecting each new vertex to exactly two of the new vertices
attached to $e_{i+1}$ via an edge contained in $T_{i+1}$.  See
Figure~\ref{fig:hyprefine}. The new vertices and edges form a binary
tree, which is VEL- and EEL-hyperbolic.  Since the refined graph
contains a hyperbolic graph, it must be hyperbolic by the
monotonicity property. The refined graph may be made into a disk
triangulation graph by adding more edges in the $T_i$ to divide the
quadrilaterals formed between the branches of the tree into
triangles.\end{proof}

\begin{theorem} \label{th:pararefine}For any infinite graph $G$ it is possible to attach vertices to the edges of $G$
to form a graph that is VEL-parabolic.\end{theorem}
\begin{proof} Let $G$ be an infinite graph and let $\{A_i\}_{i\geq 0}$ be
a collection of disjoint finite sets of edges with the property that
every path $\gamma\in\Gamma_\infty(G)$ intersects each of the $A_i$
along at least one edge.  For example, we make take $A_i$ to be the
set of edges with one endpoint in the sphere of radius $2^i$ and the
other in the sphere of radius $2^{i+1}$. Let $k_i$ be the number of
edges in $A_i$ and consider the graph $\zeta G$ that adds $k_i$
vertices incident to each edge in $A_i$ for every $i\geq 0$, and
leaves any edges not in any $A_i$ untouched. The trick is to let the
size of the $A_i$'s prescribe exactly how much to slow the vertex
growth.

Define a metric $\mu$ on $\zeta G$ by
\begin{displaymath}\mu(v) = \left\{
  \begin{array}{ll}
     \frac{1}{ik_i} & \textrm{if $v$ is incident to an edge in $A_i$} \\
     0 & \textrm{otherwise.}
  \end{array}
  \right.
\end{displaymath}
We show $\zeta G$ is parabolic by proving $\mu$ is an extremal
metric.  Let $\gamma'\in\Gamma_\infty(\zeta G)$.   By construction,
for every positive integer $i$ there are adjacent vertices $v_i,
w_i\in\gamma\cap A_i$. Thus for each $i\geq 0$, $\gamma'$ intersects
at least $k_i$ vertices incident to an edge in $A_i$. Hence,
$$L_\mu(\gamma') =\sum_{v\in \gamma} \mu(v) \geq \sum_{i=0}^\infty k_i \cdot\frac{1}{ik_i}
 =\sum_{i=0}^\infty\frac{1}{i}=\infty.$$

 On the other hand,
 $$\area(\mu)=\sum_{v\in V(\zeta G)} \mu(v)^2 =\sum_{i=0}^\infty
 \sum_{v\in A_i} \mu(v)^2 $$ $$=\sum_{i=0}^\infty k_i\cdot k_i\cdot
 \left(\frac{1}{ik_i}\right)^2 =\sum_{i=0}^\infty \frac{1}{i^2} < \infty.$$
 The two $k_i$'s beginning the second line reflect each of the $k_i$ edges
 of $A_i$ being refined $k_i$ times.
 This makes $\mu$ a finite area metric such that all paths to
 infinity have infinite length, hence $\zeta G$ has a parabolic extremal
 metric.
 \end{proof}
The construction described adds only vertices.  We now sketch how to
mimic this effect on a refined triangulations.

\begin{figure}
\begin{center}
        \includegraphics[width=3in]{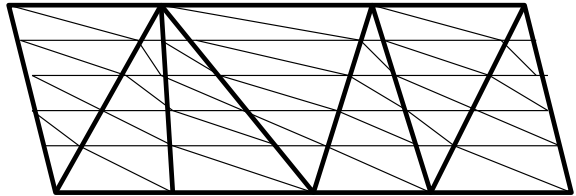}
\end{center}
\caption{The zig-zag refinement $\zeta_4$. The horizontal lines are
the levels. The diagonal lines ensure that the refinement is a
triangulation.} \label{fig:zigzag}
\end{figure}

\begin{theorem}  \label{th:pararefinetri}Let $G$ be a $VEL$-hyperbolic disk triangulation
graph.  Then there is a refinement of $G$ that yields a
$VEL$-parabolic disk triangulation graph.
\end{theorem}
\begin{proof}Let $A_i$ be the
set of edges with one endpoint in the sphere of radius $2^i$ and the
other in the sphere of radius $2^i+1$. Consider the ``zig-zag''
refinement $\zeta_n$ depicted in Figure~\ref{fig:zigzag}.  This
refinement forms a triangulation in which $n$ ``levels'' are added
between two cycles. Let $\zeta G$ be the graph formed by applying
the $\zeta_{k_i}$ refinement to the $A_i$, where $k_i$ is the number
of edges in $A_i$, and leaving untouched any edge not in some $A_i$.
 The proof now proceeds exactly
as in Theorem~\ref{th:pararefine}. Details are left to the reader,
or see \cite{woodthesis}.
\end{proof}

\section{Extensions and applications}\label{sec:dual}

\subsection{General cell complexes}

The proof of Theorem~\ref{th:vel}  depends in an essential way on
the assumption that $G$ was a triangulation.  The significant
feature of triangulations is that when a refined path cuts through a
face near some vertex $v$, it must leave again near a vertex
adjacent to $v$.  This  is no longer the case if we allow complexes
with non-triangular faces, in which refined paths may cut through
the face and emerge near vertices too far away from its entry point
to be counted properly.  In other words, our little directed graph
trick no longer works.

We can still say something if we permit some assumptions on the
planar complex $G$. For a cell complex $G=(V,E,F)$, we define its
\emph{dual complex} $G^*=(F,E^*,V)$ to be the complex whose vertices
are the faces of $G$ and whose edges are the pairs of faces of $G$
sharing an edge of $G$. We say $G$ has \emph{$(d,a)$-dually bounded
degree} if $G$ has degree $d$ and the dual complex $G^*$ has degree
$a$. The latter condition is equivalent to requiring that the faces
of $G$ each have at most $a$ sides.  For example, a bounded degree
triangulation is $(d,3)$-dually bounded for some $d$. A graph is
\emph{dually bounded} if it is $(d,a)$-dually bounded for some
$d,a$.  Our definitions for bounded refinements of cell complexes
are just as for refinements of triangulations defined in
Section~\ref{sec:refinement}.

We now offer results for dually bounded complexes similar to those
we have already obtained for triangulations.  The arguments are also
similar and we leave the reader to fill out the proof sketches below
or see \cite{woodthesis}.

 \begin{theorem}\label{th:velcomplex} Let $A=(V,E,F)$ be a $(d,a)$-dually bounded finite annular complex and
 $r$ a $b$-weakly bounded refinement of $A$.
 Then $\VEL(rA)$ and $\VEL(A)$ are $k$-comparable for some $k$
 depending only on $a$, $b$, and $d$.
 \end{theorem}
 \begin{proof}
 We assume for convenience that $b>0$ (if $b=0$, we may take
$b=1$).

Let $\mu'$ be an extremal vertex metric on $rA=(V',E',F')$. Define
$\mu(v)$ to be the greater of $\mu'(v)$ and the maximum value of
$\mu'(v')$ taken over all $v'\in V'$ incident to some edge $[v,w]$,
$w\in V$.  An argument similar to that of Theorem~\ref{th:vel} gives
$$ \VEL(rA)\leq \frac{(b+1)^2}{2}\VEL(A).$$  Conversely,
let $\mu$ be an extremal vertex metric on $A=(V,E,F)$. Let
 $rA=(V',E')$ and for each $v'\in V'$ define the \emph{face neighbors} of $v'$ as the set
 $f(v')=\{w'\in V': w',v'$ lie incident to or within the boundary of some
 face of $A\}.$
 Define a vertex metric $\mu'$ on $rG$ by
\begin{displaymath}\mu'(v') = \left\{
  \begin{array}{ll}
     \max_{w\in f(v')}(\mu(w)) & \textrm{if $v'\in V$ or $v'$} \\
     & \textrm{is incident to some edge of $E$.} \\ \\
     0 & \textrm{otherwise.}
  \end{array}
  \right.
\end{displaymath}
Again,  the argument in the proof of Theorem~\ref{th:vel} gives
$$\VEL(rA)\geq \frac{1}{a^3bd}\VEL(A),$$  Proving the theorem
with $k=\max(\frac{1}{2}(b+1)^2,a^3bd)$.\end{proof}

Degree is not a problem for edge extremal length, but we still
require a strongly bounded refinement and a bounded degree dual.

\begin{theorem} Let $A$ be a finite annular complex and
 $r$ a $(b,c)$-strongly bounded refinement of $G$.  Suppose every face of $G$ has at most $a$ sides.
 Then $\EEL(rG)$ and $\EEL(A)$ are $k$-comparable for some $k$
 depending only on $a$, $b$, and $c$.\end{theorem}
 \begin{proof}
 The proof of the relation
$$\EEL(rA)\leq
b^2\EEL(A)$$  in Theorem~\ref{th:eel} did not depend on the cells
being triangular. We take this as proved.

For the reverse relation, let $\mu$ be an edge extremal metric on
$A$. For $e'\in E'$, define $\mu'(e')$ as
\begin{displaymath}\mu'(e') = \left\{
  \begin{array}{ll}
     \mu(e) & \textrm{if there is an $e\in E$} \\
       & \textrm{such that $e'\in  \einc_r(e)$.} \\ \\
     \max(\mu(e_1), \ldots, \mu(e_n)) & \textrm{if } e' \textrm{ lies in the interior of a face}\\
                          &             \textrm{bounded by edges } e_1,\ldots,e_n \textrm{ and} \\
     & e' \textrm{ is incident-adjacent to one of the $e_i$.} \\ \\
     0 & \textrm{otherwise.}
  \end{array}
  \right.
\end{displaymath}
We leave it to the reader to adapt the proof of Theorem~\ref{th:eel}
to obtain the relation
$$\EEL(A)\leq \frac{1}{4}a^2(2ac(b+1)+b)\EEL(rA).$$  The theorem thus holds for
$k=\max(b^2,\frac{1}{4}a^2(2ac(b+1)+b) ).$ \end{proof}

By filling out a graph with annuli as before, we get the
corresponding statements on type.
 \begin{corollary}  \label{th:typecomplex}Let $G$ be an infinite planar   complex.
\begin{enumerate} \item
 If  $r_1$ is a weakly bounded refinement of $G$ and $G$ is dually bounded, then $G$ and $rG$ have the same VEL type.
\item  If  $r_2$ is a strongly bounded refinement of $G$ and $G^*$
has bounded degree, then $G$ and $rG$ have the same EEL type.

 \end{enumerate}
 \end{corollary}

Note that if a graph is dually bounded then its VEL and EEL types
are the same and we may meaningfully refer to \emph{the}
combinatorial type of the graph.

\subsection{Dual graphs}

Corollary~\ref{th:typecomplex} may be applied to relate the
combinatorial type of a complex to that of its dual.

\begin{figure}\centering

           \includegraphics[width=2in]{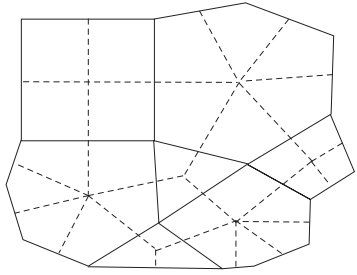}

 \caption{The graph $G^\diamond$ is a refinement of both $G$ and $G^*$.  Solid lines
 indicate edges of a complex $G$, dashed lines indicate the dual
 complex $G^*$.
 \label{fig:dual}}
\end{figure}

\begin{theorem}  A dually bounded planar cell complex has the same combinatorial type as its
 dual complex $G^*$.\end{theorem}
\begin{proof}
Let $G=(V,E,F)$ be a dually bounded planar cell complex. Consider
the refinement $G^\diamond$ constructed by adding a vertex $v_f$
inside each face of $f\in F$ and a vertex $v_e$ incident to each
edge $e\in E$.  Connect these new vertices by adding edges of the
form $[v_f, v_e]$, where $e$ is an edge bounding the face $f$.
Roughly, we are superimposing $G^*$ onto $G$ and attaching vertices
where they intersect.  See Figure~\ref{fig:dual}. This construction
clearly defines a strongly bounded refinement $rG=G^\diamond$ of
$G$, so $G^\diamond$ shares its type.  We associate $v_f$ to its
dual vertex $f\in V^*=F$, a pair of edges $[v_f, v_e],[v_g, v_e]$,
$e\in E$, $f,g\in F$ to the dual edge $[f,g]\in E^*$, and we note
that each vertex $v\in V$ lies inside a distinct face  formed by the
edge pairs associated to $E^*$. With these associations, we see that
$G^\diamond$ is isomorphic not only to $rG$, but also to $r(G^*)$ --
the vertices $v_e$ attached to an edge $e$ are identified with
vertices $v_{e^*}$ attached to the dual edges, and we similarly
reverse the roles of vertices and faces.  We have
$G^\diamond=rG\cong r(G^*)$, and so $G$ and $G^*$ have the same type
by Corollary~\ref{th:typecomplex}.\end{proof}

\subsection{Outer Spheres}

The refinement theorems we have developed suggest solving  the
combinatorial type problem for a specific graph $G$ by finding a
dually bounded subcomplex of $G$ from which $G$ can be obtained by a
bounded refinement.  The purpose of this section is to construct a
candidate subgraph.

Let $G=(V,E)$ be a disk triangulation graph with distinguished base
vertex $v_0$.  Define the \emph{outer sphere of radius $n$}
$\osph(n)$ to be the collection of vertices $v\in G$ such that
$|v|=n$ and for which there is a path $\gamma^+(v)$ from $v$ to
infinity containing no other vertices of norm $n$.  The edges of
$\osphn$ are the edges of $E$ whose vertices both lie in $\osphn$.

The typical spheres $S(n)=\{v:|v|= n\}$ may be massively
disconnected because of geodesics that cannot be extended, mimicking
a classical phenomenon. The outer spheres ignore these unwanted
components. We  show that they are cycle graphs.

\begin{lemma} $\osph(n)$ is a cycle graph for all
integers $n>0$. \end{lemma}

\begin{proof} Fix $n$ and let $v\in\osph(n)$. Note that $v$ is the only
vertex in $\gamma^+(v)$ whose norm is not strictly larger than $n$.
To see this, observe that since $\gamma^+(v)$ goes to infinity, it
must contain vertices of arbitrarily large norm. Were there a vertex
of norm less than $n$, then the fact that adjacent vertices differ
in norm by at most 1 implies that there would also be a vertex of
norm $n$, contradicting the definition of $\gamma^+(v)$.

Similarly, for every vertex $v\in\osph(n)$ there is a path
$\gamma^-(v)$ connecting $v_0$ to $v$ such that $\gamma^-(v)$ is of
length $n$ and therefore contains no other vertices with norm
greater than or equal to $n$.  This is immediate from the
definitions of $\osph(n)$ and norm. Altogether, the concatenated
path $\gamma(v)= \gamma^-(v)\cup\gamma^+(v)$ connects the base point
$v_0$ to infinity so that $v$ is the only vertex in the path of norm
$n$, all vertices before $v$ in $\gamma(v)$ have norm less than $n$,
and all vertices after $v$ have norm larger than $n$.

\begin{figure}
\center
\includegraphics[width=2in]{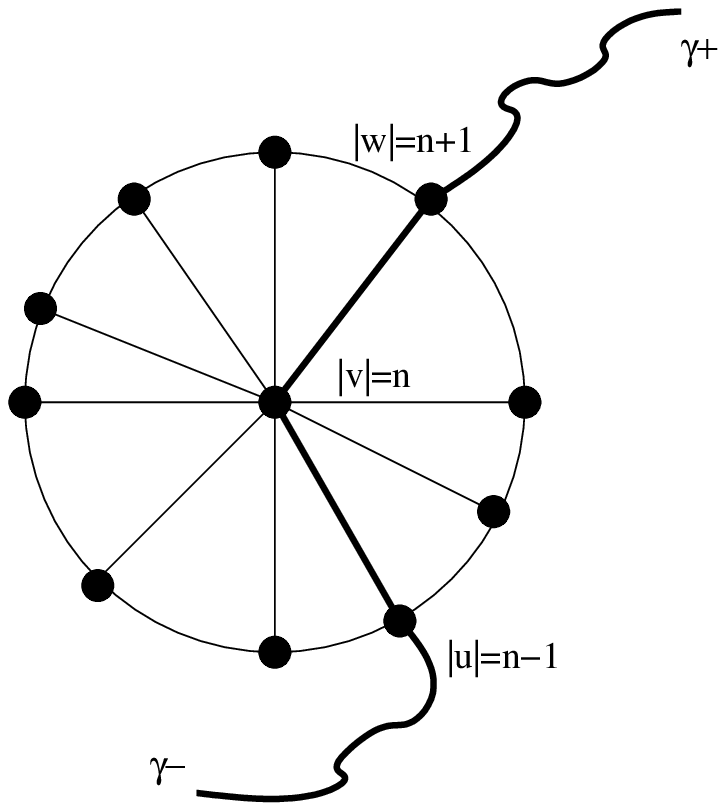} \caption{A vertex in $\osphn$ has two neighbors in $\osphn$}
\label{fig:deggt1}
\end{figure}

The successor $w$ and predecessor $u$ of $v$ in $\gamma(v)$ must
therefore have norms $n+1$ and $n-1$, respectively.  Since $G$ is a
triangulation, the set of neighbors of $v$ are cyclicly connected by
edges and so the closed loop connecting the neighbors of $v$
contains at least one vertex each of norm $n+1$ and $n-1$. These
vertices divide the loop of neighbors of $v$ into two segments, each
containing a vertex with norm $n$, again because successive vertices
along a path may differ in norm by at most 1. So $v$ has at least
two neighbors $v_1, v_2$ with norm $n$.  Start at $v_1$ and proceed
around the loop of edges toward $w$. We may assume $v_1$ was chosen
so that no other vertices of norm $n$ are encountered.  Then travel
along $\gamma^+(v)$ away from $v$, giving a path from $v_1$ to
infinity that contains no vertices with norm $n$.  Repeating for
$v_2$, we have shown that any $v$ in $\osphn$ has at least two
neighbors in $\osphn$, i.e. that every vertex in $\osphn$ has degree
at least two. See Figure~\ref{fig:deggt1}.

Now suppose for contradiction that $v\in\osphn$ has three neighbors
$w_1,w_2,w_3$ in $\osphn$ connected to $v$ by edges $e_1,e_2,e_3$.
Then to each $w_i$ there is a path $\gamma_i^+$ connecting $w_i$ to
infinity whose interior vertices all have norm at least $n$. Assume
without loss of generality that the $\gamma_i^+$ do not contain the
base point $v_0$.  The set
$\cc\setminus(\gamma_1^+\cup\gamma_2^+\cup\gamma_3^+\cup\{e_1,e_2,e_3\})$
is a collection of regions in the plane whose boundaries contain
only vertices of norm greater than or equal to $n$ and at most two
of the $w_i$. See Figure~\ref{fig:deglt3}.

\begin{figure}
\center
\includegraphics[height=1in,width = 2in]{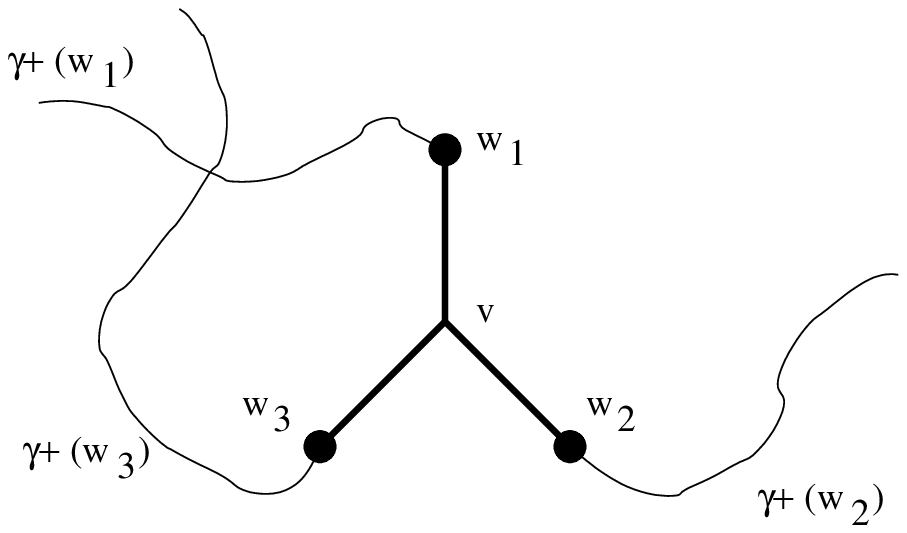} \caption{A vertex in $\osphn$ cannot have three neighbors in $\osphn$}
\label{fig:deglt3}
\end{figure}

Consider the region $R$ containing the base point $v_0$ and suppose
$w_1$ is not in its boundary $\partial R$.  Then by assumption there
is a path $\gamma^-(w_1)$ connecting $v_0$ to $w_1$ and containing
only vertices with norm less than $n$ except the endpoints. But
$\partial R$ separates $v_0$ from $w_1$, so $\gamma^-(w_1)$ must
intersect $\partial R$, all of whose vertices have norm greater than
or equal to $n$. This is a contradiction, so $v$ cannot have degree
greater than two.  Since $\osphn$ is compact and all of its vertices
have degree two, $\osphn$ must be a union of disjoint cycle graphs.
We have only to show $\osphn$ has but one component.

We begin by showing that if $S$ is a component of $\osphn$, then
$v_0$ is contained in the bounded component of $\cc\setminus S$.
Suppose for contradiction that $v_0$ is in the unbounded component
and let $v\in S$.  As before, there is a path $\gamma^-(v)$ from
$v_0$ to $v$ with all interior vertices having norm smaller than
$n$, and a path $\gamma^+(v)$ from $v_0$ to infinity with all
interior vertices of norm larger than $n$.  Then the paths
$\gamma^+(v)$ and $\gamma^-(v)$ both lie entirely entirely in the
unbounded component except where they meet at $v$. But then the
successor of $v$ in $\gamma^+(v)$ and the predecessor of $v$ in
$\gamma^-(v)$ both lie in the unbounded component of $\cc\setminus
S$. These vertices have norm $n+1$ and $n-1$, respectively, and by
the argument used above, we have that there must be another element
of $\osphn$ lying in the unbounded component of $S$, a contradiction
to the assumption that $S$ is a cycle. See Figure~\ref{fig:oscctd}.

\begin{figure}
\center
\includegraphics[width=2in]{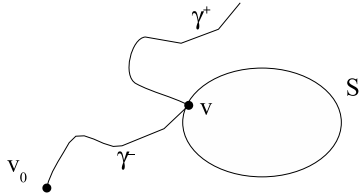} \caption[$\osphn$ is connected]{Both paths lie in the unbounded component of $\cc\setminus S$}
\label{fig:oscctd}
\end{figure}

The only remaining possibility is that the components of $\osphn$
are concentric.  But if there is more than one component, then it is
not possible to find a path to infinity from a vertex in an inner
component without crossing the outermost component which contains
only vertices of norm $n$.  We are left to conclude that $\osphn$ is
connected, proving the claim.
\end{proof}

For a disk triangulation graph $G$, its outer spheres suggest a
subgraph of $G$ for study.  Define the \emph{outer sphere skeleton}
$G_\mathcal{O}$ to be the union of the outer spheres $\osphn$ along
with all edges of the form $[v_n, v_{n+1}]$ where $v_n\in \osphn$
and and $v_{n+1}\in\osph(n+1)$.

$G_\mathcal{O}$ is simply the set of outer spheres along with the
edge geodesics connecting the outer spheres to $v_0$. It discards
the isolated face subdivisions that we have already seen cannot
impact VEL type for dually bounded degree complexes. $G_\mathcal{O}$
has an appealing structure. All vertices on $\osphn$ may be traced
back to the base vertex $v_0$ by working backward through each of
the previous outer spheres.
 Each face
 of $G_\mathcal{O}$ lies between two outer spheres $\osphn$ and $\osph(n+1)$.  The face
 is bounded by two edges connecting these spheres, at most one
 edge of $\osph(n+1)$, and any number of consecutive edges along
 $\osphn$.  The following is a special case of
 Theorem~\ref{th:typecomplex}.
\begin{theorem}
Let $G$ be a disk triangulation graph.  If the outer sphere skeleton
$G_\mathcal{O}$ is VEL-hyperbolic, then $G$ is VEL-hyperbolic.  If
$G_\mathcal{O}$ is dually bounded and VEL-parabolic, then $G$ is
VEL-parabolic.\end{theorem}

\section{Growth and type\label{sec:growth}}

We now consider to what extent combinatorial type can be determined
from the growth of a graph.  There are several results that say
roughly that slow-growing graphs are parabolic and fast-growing
graphs are hyperbolic.  For example, Rodin and Sullivan \cite{rs}
showed using circle packings (which connect to vertex extremal
length by \cite{heschrammhp}) that linear growth of spheres in a
bounded degree graph implies parabolicity, whereas Siders
\cite{siders} used electrical methods (via \cite{duffin}) to
determine the type of a graph formed by interspersing cycles of 6-
and 7-degree vertices.  See also \cite{woess}. In this section, we
establish a sharp parabolicity condition that requires no symmetry
or degree restrictions of the graph.

Define $\log_{(0)}(x)=x$ and $\log_{(m+1)}(x)=\log(\log_{(m)}(x))$.
Inductive application of the chain rule gives the derivative
$\log'_{(m+1)}(x)=\frac{d}{dx}\log_{(m+1)}(x)=(\log_{(0)}(x)
\log_{(1)}(x)\cdots\log_{(m)}(x))^{-1}$, and define the
\emph{generalized $p$-series} $$\wp(m,p) = \sum_{j=e_m}^\infty
\frac{\log'_{(m)}(j)}{(\log_{(m)}(j))^p}$$
$$ = \sum_{j=e_m}^\infty \frac{1}{j(\log j)( \log\log j)\cdots
(\log_{(m-1)}j )(\log_{(m)}j)^p}$$ where $e_m$ is the smallest
integer for which $\log_{(m)}(e_m)$ is defined.

\begin{lemma} \label{th:pseries} Let $m$ be a nonnegative integer and $c>0$ such that $\log_{(m)}(c)>0$. \begin{enumerate} \item The integral $$\int_c^\infty
\frac{\log'_{(m)}(x)}{(\log_{(m)}(x))^p}dx$$ converges if and only
if $p>1$. \item The generalized $p$-series $\wp(m,p)$ converges if
and only if $p>1$.\end{enumerate}\end{lemma}
\begin{proof}
Making the substitution $u=\log_{(m)}(x)$ and $du=\log'_{(m)}(x)\
dx$, we integrate
$$\int_c^\infty
\frac{\log'_{(m)}(x)}{(\log_{(m)}(x))^p}dx=\int_{\log_{(m)}(c)}^\infty
\frac{1}{u^p}du,$$ which converges if and only if $p>1$.  This
proves the first part of the lemma, and the second part follows
because the generalized $p$-series approximates this integral. See
\cite{rudin} for an alternate proof and general
discussion.\end{proof}

We can now establish the parabolicity condition.

\begin{theorem}\label{th:para}
Let $G=(V,E)$ be a graph with distinguished vertex $v_0$.  Suppose
there is  a collection of subgraphs $\{C_j=(V_j,E_j)\}$ of $G$ such
that for each $j\in\mathbb{N}$, every vertex path from $v_0$ to
infinity intersects $C_j$ at least once.  If there is a positive
integer $m\in\mathbb{N}$ and  a positive constant $K>0$ such that
for every $j\geq m$ we have $|V_j|\leq K j(\log j)( \log\log
j)\cdots (\log_{(m-1)}j)=\frac{K}{\log'_{(m)}j},$ then $G$ is
VEL-parabolic.
\end{theorem}

\begin{proof}
 Define a vertex metric $\mu$ on $G$ by
\begin{displaymath}\mu(v) = \left\{
  \begin{array}{ll}
     \frac{\log'_{(m)}(j)}{\log_{(m)}(j)} & \textrm{if $v\in V_j$ and $j\ge e_m$} \\
     0 & \textrm{otherwise}
  \end{array}
  \right.
\end{displaymath}
 Let $\Gamma$ be the collection
of transient vertex paths in $G$ based at $v_0$. Since each
$\gamma\in\Gamma$ must contain at least one vertex in each of the
$V_j,$  we have
$$L_\mu(\gamma)=\sum_{v\in\gamma}\mu(v)\geq \sum_{j=e_m}^\infty \frac{\log'_{(m)}(j)}{\log_{(m)}(j)}=\wp(m,1).$$ This
sum  diverges by Lemma~\ref{th:pseries}, showing all paths in
$\Gamma$ have infinite length.

We now show that $\mu$ has finite area.
$$\area(\mu)=\sum_{v\in V}\mu(v)^2 = \sum_{j=e_m}^\infty
\frac{|V_j|(\log'_{(m)}(j))^2}{(\log_{(m)}(j))^2}$$ $$ \leq
 \sum_{j=e_m}^\infty\frac{K\log'_{(m)}(j)}{\log_{(m)}(j))^2}
=K\wp(m,2)$$ which is finite by Lemma~\ref{th:pseries}.

It follows that $\mu$ is a parabolic extremal vertex metric for $G$
and so $G$ is VEL-parabolic.

\end{proof}

Of particular note in Theorem~\ref{th:para} is the case where $G$ is
a disk triangulation graph and  $C_j$ is the outer sphere of radius
$j$ about $v_0$. Note also that the theorem does not hold in the EEL
context as stated, but the proof can be adapted easily if we require
slow growth of the number of edges emanating from the cycles $C_j$.
Details are an exercise.

There is no hope for a converse to Theorem~\ref{th:para}, as P.
Soardi \cite{soardi} has a parabolic graph with exponential growth.
The example has bounded degree and thus functions in both the VEL
and EEL settings.

We can get a sharpness result, however.

\begin{theorem}  Let $\{a_i\}_{i=0}^\infty$ be a sequence of positive integers
such that $\sum \frac{1}{a_i} < \infty$.  Then there is a
VEL-hyperbolic disk triangulation graph with a vertex $v_0\in G$
such that the $n$-sphere about $v_0$ is a cycle containing at most
$a_n$ vertices.
\end{theorem}

\begin{figure}\centering

 \subfigure[\textsf{A square tiling with $n^2$ squares in the $n$th row}]
 {\hspace{1.0in}\includegraphics[width=2in]{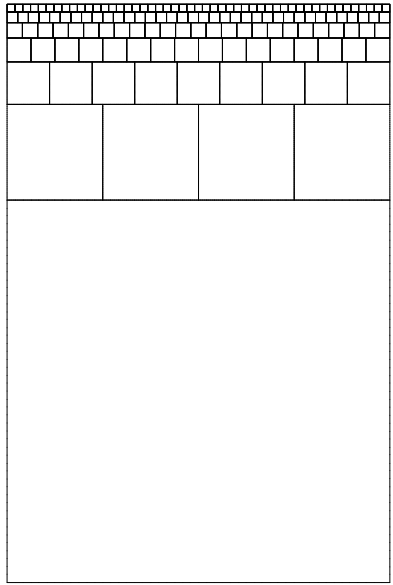}\hspace{1.0in}\label{fig:quad1}}

 \subfigure[\textsf{The corresponding triangulation}]
 {\hspace{1.0in}\includegraphics[width=2in]{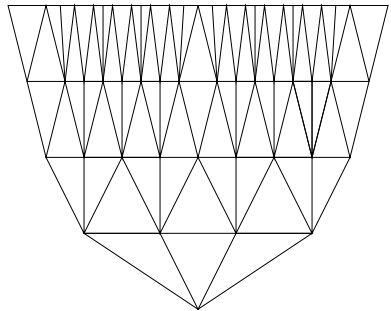}\hspace{1.0in}\label{fig:quad2}}

\caption[Quadratic growth in a hyperbolic triangulation]{A
hyperbolic triangulation with quadratic growth \label{fig:quad}}
\end{figure}

\begin{proof} Assume for convenience that $\{a_i\}$ is monotone.
The construction uses a connection of vertex extremal length to
square tilings established in \cite{cfp} and \cite{schrammsquare}:
if $G=(V,E)$ is a triangulation of a finite closed quadrilateral,
then there is a tiling by squares of a rectangle so that each vertex
in $V$ corresponds to a square in the tiling and two squares
intersect if and only if their corresponding vertices determine an
edge in $E$.  The side lengths of the squares determine an extremal
metric for $G$, and the vertex extremal length of $G$ is the aspect
ratio of the rectangle.

Construct a tiling $T_n, n>1$ as follows.  Begin with a unit square,
which  corresponds to $v_0$.  Add a row of $a_1$ squares with side
lengths all $\frac{1}{a_1}$ along one side of the unit square.
Continuing adding rows of $a_i$ squares with side lengths
$\frac{1}{a_i}$ for all $i\leq n$.  Figure~\ref{fig:quad}
illustrates the construction for $a_n = (n+1)^2$.  Minor adjustments
to some of the $a_i$ may be necessary in the unlucky case that some
$a_i$ divides $a_{i+1}$ and the corresponding graph is not a
triangulation.

The sides of the tiled rectangle are $1$ and $1+\sum_{i=1}^n
\frac{1}{a_i}$, so the vertex extremal length of the corresponding
triangulation is $1+\sum \frac{1}{a_i}.$  By assumption, this
remains finite as $n\to \infty$ and so the limiting triangulation
$T_\infty$ is VEL-hyperbolic.  This triangulation forms the required
disk triangulation graph if we identify the pairs of boundary
vertices that lie in common spheres about $v_0$.
\end{proof}

We have thus shown, for example, that $O(n)$ sphere growth implies a
graph is parabolic, whereas $O(n^{1+\varepsilon})$ sphere growth is
indeterminate for $\varepsilon>0$.

Note that the construction also works in the EEL setting because
VEL-hyperbolic implies EEL-hyperbolic.

\bibliographystyle{amsalpha}
\bibliography{typepaper}

\end{document}